\documentclass[a4paper, 11pt]{amsart}

\usepackage[english]{babel}
\usepackage[utf8x]{inputenc}
\usepackage[T1]{fontenc}

\usepackage[a4paper,top=4cm,bottom=4cm,left=3cm,right=3cm,marginparwidth=1.75cm]{geometry}

\usepackage{amsmath}
\usepackage{amssymb}
\usepackage{graphicx}
\usepackage{mathrsfs}
\usepackage[colorinlistoftodos]{todonotes}
\usepackage[colorlinks=true, allcolors=blue]{hyperref}
\usepackage[tableposition=top]{caption}
\DeclareCaptionLabelFormat{gostfigure}{Fig #2}
\DeclareCaptionLabelSeparator{gost}{~---~}
\newtheorem{theorem}{Theorem}[section]
\newtheorem{definition}{Definition}[section]
\newtheorem{proposition}{Proposition}[section]
\newtheorem{corollary}{Corollary}[section]
\newtheorem{hypothesis}{Hypothesis}[section]
\title{On a category of cotangent sums related to the Nyman-Beurling criterion
for the Riemann Hypothesis}
\date{}
\author{Derevyanko Nikita}
\address{Moscow Institute of Physics and Technology,
141700 Dolgoprudny,
Institutskiy per, d. 9,
Russia.
}
\email{nikita.derevyanko@phystech.edu}
\author{Kovalenko Kirill}
\address{National Research University Higher School of Economics, 101000 Moscow, Myasnitskaya ulitsa, d.20,  Russia
}
\email{kdkovalenko@edu.hse.ru}

\begin{document}
\maketitle

\begin{abstract}
  The purpose of the present paper is to provide a general overview of a variety
  of results related to a category of cotangent sums which have been proven to
  be associated to the so-called Nyman-Beurling criterion for the Riemann
  Hypothesis. These sums are also related to the Estermann Zeta function.
\end{abstract}

\section{Introduction}
This paper is focused on applications of certain cotangent sums to different problems
related to the Riemann Hypothesis. The expression for the sums in question is the following

\bigskip
\begin{definition}
  \begin{equation}\label{eq:c0}
    c_0\left(\frac{r}{b}\right) := - \sum_{m=1}^{b-1}\frac{m}{b}\cot\left(\frac{\pi mr}{b}\right),
  \end{equation}
  where $r,b \in \mathbb{N},\ b \ge 2,\ 1 \le r \le b$ and $(r,b) = 1$.
\end{definition}

\subsection{Nyman–Beurling criterion for the Riemann Hypothesis}

There are many interesting results concerning these cotangent sums, but
initially we will present some general information about the Riemann Hypothesis
and some related problems. Moreover, our aim is to provide motivation for the use of
cotangent sums in these problems.

In this paper we shall denote a complex variable by $s = \sigma + it$, where
$\sigma$ and $t$ are the real and imaginary part of $s$ respectively.

\bigskip
\begin{definition}
  The Riemann zeta function is a function of the complex variable $s$ defined in
  the half-plane $\{\sigma > 1\}$ by the absolutely convergent series
  \begin{equation}
    \zeta(s) := \sum_{n=1}^{\infty} \frac{1}{n^s}.
  \end{equation}
\end{definition}

As shown by B. Riemann, $\zeta(s)$ extends to $\mathbb{C}$ as a meromorphic function
with only a simple pole at $s=1$, with the residue 1, and satisfies the
functional equation
\begin{equation}
  \zeta(s) = 2^s\pi^{s-1}\sin\left(\frac{\pi s}{2}\right)\Gamma(1-s)\zeta(1-s).
\end{equation}

For negative integers, one has a convenient representation of the Riemann zeta function
in terms of Bernoulli numbers:
$$
  \zeta(-n) = (-1)^n\frac{B_{n+1}}{n+1},\ \text{for}\ n \geq 0.
$$
By the above formula one can easily deduce that $\zeta(s)$ vanishes when $s$ is a negative
even integer because $B_m = 0$ for all odd $m$ other than $1$. The negative even integers 
are called trivial zeros of the Riemann zeta function. All other complex points 
where $\zeta(s)$ vanishes are called non-trivial zeros of the Riemann
zeta function, and they play a significant role in the distribution of primes.

The actual connection with the distribution of prime numbers was observed in Riemann's
1859 paper. It is in this paper that Riemann proposed his well known hypothesis.
\begin{hypothesis}[Riemann]
    The Riemann zeta function $\zeta(s)$ attains its non-trivial zeros only in complex points with $\sigma = \frac{1}{2}$. The line on the complex plane given by the equation $\sigma = \frac{1}{2}$ is usually called "critical".
\end{hypothesis}

The Nyman–Beurling–Baez–Duarte–Vasyunin (also simply known as Nyman–Beurling)  approach to the Riemann Hypothesis enables us
to associate the study of cotangent sums and the Riemann Hypothesis
through the following theorem:

\begin{theorem}\label{th:NB}
  The Riemann Hypothesis is true if and only if
  $$
    \lim_{N\rightarrow+\infty} d_N = 0,
  $$
  where
  \begin{equation}\label{eq:dN}
    d_N^2 = \inf_{D_N}\frac{1}{2\pi}\int_{-\infty}^{+\infty}\left| 1-
    \zeta\left(\frac{1}{2} + it\right)D_N\left(\frac{1}{2} + it\right)\right|^2 \frac{dt}{\frac{1}{4}+t^2}
  \end{equation}
  and the infimum is taken over all Dirichlet polynomials
  \begin{equation}\label{eq:D_N}
    D_N(s) = \sum_{n=1}^N \frac{a_n}{n^s},\ a_n \in \mathbb{C}.
  \end{equation}
\end{theorem}

In his paper \cite{Bagchi}, B. Bagchi used a slightly different formulation of
Theorem \ref{th:NB}. In order to state it, we have to introduce some definitions.

\bigskip
\begin{definition}
  The Hardy space $H^2(\Omega)$ is the Hilbert space of all analytic functions $F$ on
  the half-plane $\Omega$ (we define it for a right half-plane $\{\sigma > \sigma_0\}$ of the 
  complex plane) such that
  $$
    \|F\|^2 := \sup_{\sigma > \sigma_0}\frac{1}{2\pi}
    \int_{-\infty}^{+\infty}|F(\sigma + it)|^2dt < \infty.
  $$
  Everywhere in this section we will use $\Omega = \{\sigma > \frac{1}{2}\}$.
\end{definition}

\bigskip
\begin{definition}
  For $0 \le \lambda \le 1$, let $F_{\lambda} \in H^2(\Omega)$ be defined by
  $$
    F_{\lambda}(s) = (\lambda^s - \lambda)\frac{\zeta(s)}{s}, s \in \Omega,
  $$
  and for $l = 1, 2, 3, \ldots$, let $G_l \in H^2(\Omega)$ be defined by
  $G_l = F_{\frac{1}{l}}$, i.e
  $$
    G_l(s) = (l^{-s} - l^{-1})\frac{\zeta(s)}{s}, s \in \Omega.
  $$
  Also, let $E \in H^2(\Omega)$ be defined by
  $$
    E(s) = \frac{1}{s}, s \in \Omega.
  $$
\end{definition}

Now we can state the reformulation of Theorem \ref{th:NB} which was used in paper
\cite{Bagchi}.

\bigskip
\begin{theorem}
  The following statements are equivalent:
  \begin{enumerate}
    \item The Riemann Hypothesis is true;
    \item $E$ belongs to the closed linear span of the set
    $\{G_l: l=1,2,3 \ldots\}$;
    \item $E$ belongs to the closed linear span of the set
    $\{F_{\lambda}: 0 \le \lambda \le 1\}$.
  \end{enumerate}
\end{theorem}

The plan of the proof is to verify three implications:
$1 \rightarrow 2$, $2 \rightarrow 3$ and $3 \rightarrow 1$.
\bigskip

The first implication is the most challenging of all three. It is proven using
some famous results obtained under the assumption that the Riemann Hypothesis is true, among which are
Littlewood's theorem \ref{th:Litt} and the Lindel\"{o}f hypothesis \ref{def:Lind},
and some standard techniques of functional analysis, particularly concerning
convergence in the norm. More details can be found in the
original paper by B. Bagchi \cite{Bagchi}.

\bigskip
\begin{hypothesis}[Lindel\"{o}f]\label{def:Lind}
  If the Riemann Hypothesis is true, then
  $$
    \forall \varepsilon > 0,\ \zeta\left(\frac{1}{2} + it\right) = O(t^{\varepsilon}).
  $$
\end{hypothesis}

\bigskip

\textbf{Remark:} A very interesting and novel approach to the Lindel\"{o}f  hypothesis is presented by A. Fokas \cite{Fokas}.

\bigskip

\begin{theorem}[Littlewood]\label{th:Litt}
  If the following conditions are satisfied:
  \begin{itemize}
    \item $
      \lim_{r \rightarrow 1^{-}}\sum_{n=0}^{\infty}r^nc_n = a,
      \ \forall i\ c_i \in \mathbb{C},\ a \in \mathbb{C},
    $
    \item $
      c_n = O(\frac{1}{n}),
    $
  \end{itemize}

  then
  $$
    \sum_{n=0}^{\infty}c_n = a.
  $$
\end{theorem}

The second implication follows from the embedding
$$
\{G_l: l=1,2,3 \ldots\} \subset \{F_{\lambda}: 0 \le \lambda \le 1\}.$$

\bigskip

To prove the third implication ($3 \rightarrow 1$), suppose that the Riemann Hypothesis is false. Then
$\exists s_0 = \sigma_0 + it_0: \zeta(s_0) = 0$ and $\sigma_0 \neq \frac{1}{2}$,
which implies that $\forall 0 \le \lambda \le 1\ F_{\lambda}(s_0) = 0$. That
together with statement 3 gives that
$E(s_0) = \frac{1}{s_0} = 0$, i.e $0 = 1$. This contradiction
completes the proof.

\subsection{The cotangent sum's applications to problems related to the Riemann Hypothesis}

The main motivation behind the study of the cotangent sum \eqref{eq:c0} follows
from Theorem \ref{th:NB}, which constitutes an equivalent form of the Riemann Hypothesis.

Asymptotics for $d_N$ \eqref{eq:dN} under the assumption that the Riemann Hypothesis is true have been
studied in several papers. S. Bettin, J. Conrey and D. Farmer in
\cite{Bettin_Conrey_Farmer} obtained the following result.

\begin{theorem}
  If the Riemann Hypothesis is true and if
  $$
    \sum_{|Im(\rho)|\leq T} \frac{1}{|\zeta^{\prime}(\rho)|^2}
    \ll T^{\frac{3}{2} - \delta}
  $$
  for some $\delta > 0$, with the sum on the left hand side taken over all distinct zeros $\rho$ of the Riemann zeta function with imaginery part less than or equal to $T$, then
  $$
    \frac{1}{2\pi} \int_{-\infty}^{+\infty}\left| 1-
    \zeta\left(\frac{1}{2} + it\right)V_N\left(\frac{1}{2} + it\right)\right|^2
    \frac{dt}{\frac{1}{4}+t^2} \sim \frac{2 + \gamma - \log 4\pi}{\log N}
  $$
  for
  \begin{equation}\label{eq:V_N}
    V_N(s) := \sum_{n=1}^N\left(
    1 - \frac{\log n}{\log N}\right)\frac{\mu(n)}{n^s}.
  \end{equation}
\end{theorem}
We should mention that in the sequel $\gamma$ stands for the Euler–Mascheroni constant. Also, here $\mu$ is the M\"{o}bius function.

Also, from results of \cite{Bettin_Conrey_Farmer} it follows that under some restrictions,
the infimum from \eqref{eq:dN} is attained for $D_N = V_N$.

Nevertheless, it is interesting to obtain an unconditional estimate for $d_N$.

In order to proceed further, we shall study equation \eqref{eq:dN} in more detail.
In particular, we can expand the square in the integral:

\begin{align*}
  d_N = & \inf_{D_N}\Big(
  \int_{-\infty}^{+\infty}\left(
    1-\zeta\left(\frac{1}{2} + it\right)
    D_N\left(\frac{1}{2} + it\right) -
    \zeta\left(\frac{1}{2} - it\right)
    \bar{D_N}\left(\frac{1}{2} + it\right)
    \right) \frac{dt}{\frac{1}{4}+t^2} \\
    & + \int_{-\infty}^{+\infty}
    \left|\zeta\left(\frac{1}{2} + it\right)\right|^2
    \left|D_N\left(\frac{1}{2} + it\right)\right|^2
    \frac{dt}{\frac{1}{4}+t^2}\Big).
\end{align*}

The integral in the second summand can be expressed as

\begin{equation}\label{eq:exp_integral}
  \sum_{1\leq r,b \leq N} a_r \bar{a_b} r^{-\frac{1}{2}} b^{-\frac{1}{2}}
  \int_{-\infty}^{+\infty}\left|\zeta\left(\frac{1}{2} + it\right)\right|^2
  \left(\frac{r}{b}\right)^{it} \frac{dt}{\frac{1}{4}+t^2},
\end{equation}
where $a_i$ are from the definition \eqref{eq:D_N} of $D_N$.

Therefore, the integral 
$$\int_{-\infty}^{\infty}\left|\zeta(\frac{1}{2}+it)\right|^2
\left(\frac{r}{b}\right)^{it}\frac{dt}{\frac{1}{4} + t^2}$$ 
plays an important
role in the Nyman-Beurling criterion for the Riemann Hypothesis. Moreover, one can prove
that this integral can be expressed via the so-called Vasyunin sum.

\begin{definition}
  The Vasyunin sum is defined as follows:
  \begin{equation}
    V\left(\frac{r}{b}\right) := \sum_{m=1}^{b-1}\left\{\frac{mr}{b}\right\}\cot\left(\frac{\pi mr}{b}\right),
  \end{equation}
  where $\{x\} = x - \lfloor x\rfloor, x\in\mathbb{R}$.
\end{definition}

The following proposition holds true:

\begin{proposition}\label{prop:int_to_vas}
  \begin{eqnarray}
    \lefteqn{
    \frac{1}{2\pi(rb)^{1/2}}\int_{-\infty}^{\infty}\left|\zeta\left(\frac{1}{2}+it\right)\right|^2
    \left(\frac{r}{b}\right)^{it}\frac{dt}{\frac{1}{4} + t^2}} \nonumber\\
    & & = \frac{\log2\pi - \gamma}{2}\left(\frac{1}{r} + \frac{1}{b}\right) +
    \frac{b-r}{2rb}\log\frac{r}{b} -
    \frac{\pi}{2rb}\left(V\left(\frac{r}{b}\right) + V\left(\frac{b}{r}\right)\right).
  \end{eqnarray}
\end{proposition}

One can note that the only non-explicit function on the right hand side
of this formula is the Vasyunin sum.

The next equation connects this result with the cotangent sums in question.

\begin{proposition} It holds that
  \begin{equation}
    V\left(\frac{r}{b}\right) = -c_0\left(\frac{\bar{r}}{b}\right),
  \end{equation}
  where $\bar{r}$ is such that $\bar{r}r \equiv 1 (mod\ b)$.
\end{proposition}

The cotangent sum $c_0$ can also be used to describe some special values of the
Estermann zeta function.

\begin{definition}
  The Estermann zeta function $E(s, \frac{r}{b}, \alpha)$ is defined by the Dirichlet
  series
  \begin{equation}
    E\left(s, \frac{r}{b}, \alpha\right)=
    \sum_{n \ge 1}\frac{\sigma_{\alpha}(n)\exp(\frac{2\pi inr}{b})}{n^s},
  \end{equation}
  where $Re\ s > Re\ \alpha + 1,\ b \ge 1,\ (r, b) = 1$, and
  \begin{equation}
    \sigma_{\alpha}(n)=\sum_{d|n}d^{\alpha}.
  \end{equation}
\end{definition}

One can show that the Estermann zeta function $E(s, \frac{r}{b}, \alpha)$
satisfies the following functional equation:

\begin{eqnarray}
  \lefteqn{
  E(s, \frac{r}{b}, \alpha) =
  \frac{1}{\pi}\left(\frac{b}{2\pi}\right)^{1 + \alpha - 2s}
  \Gamma(1-s)\Gamma(1 + \alpha - s)} \nonumber\\
  & & \times \left(\cos\left(\frac{\pi\alpha}{2}\right)
  E\left(1+\alpha-s, \frac{\bar{r}}{b}, \alpha\right) -
  \cos\left(\pi s - \frac{\pi\alpha}{2}\right)
  E\left(1 + \alpha-s, -\frac{\bar{r}}{b}, \alpha\right)\right),
\end{eqnarray}
where $r$ is such that $\bar{r}r \equiv 1 (mod\ b)$.

Properties of $E(0, \frac{r}{b}, 0)$ were used by R. Balasubramanian,
J. Conrey, and D. Heath-Brown \cite{BCHB} to prove an asymptotic formula for
\begin{equation}
  I = \int_{0}^{T}\left|\zeta\left(\frac{1}{2} + it\right)\right|^2\left|A\left(\frac{1}{2} + it\right)\right|^2dt,
\end{equation}
where $A(s)$ is a Dirichlet polynomial.

Asymptotic results for the function $I$ as well as other functions of this type, are useful for estimating a lower
bound for the portion of zeros of the Riemann zeta function $\zeta(s)$ on the
critical line.

The following result of M. Ishibashi from \cite{Ishibashi} concerning
$E(s, \frac{r}{b}, \alpha)$ for $s = 0$, 
provides the connection of the Estermann zeta function with the cotangent sum
$c_0(\frac{r}{b})$, simply by setting $\alpha = 0$.

\begin{theorem}[Ishibashi]
  Let $b \ge 2,\ 1 \le r \le b,\ (r,b)=1,\ \alpha \in \mathbb{N}\cup\{0\}$. Then

  for even $\alpha$, it holds that
  \begin{equation}
    E\left(0, \frac{r}{b}, \alpha\right) = \left(-\frac{i}{2}\right)^{\alpha+1}
    \sum_{m=1}^{b-1}\frac{m}{b}\cot^{(\alpha)}\left(\frac{\pi mr}{b}\right)
    + \frac{1}{4}\delta_{\alpha,0},
  \end{equation}
  where $\delta_{\alpha,0}$ is the Kronecker delta function.

\bigskip
For odd $\alpha$, it holds that
\begin{equation}
  E\left(0, \frac{r}{b}, \alpha\right) = \frac{B_{\alpha + 1}}{2(\alpha + 1)}.
\end{equation}

In the special case when $r = b = 1$, we have
$$
  E(0, 1, \alpha) = \frac{(-1)^{\alpha + 1}B_{\alpha + 1}}{2(\alpha + 1)},
$$

where $B_m$ is the $m$-th Bernoulli number and $B_{2m + 1} = 0$,
$$
  B_{2m} = \frac{(-1)^{m+1}2(2m)!}{(2\pi)^{2m}}\zeta(2m),\ for\ m \ge 1.
$$
\end{theorem}

Thus for  $b \ge 2,\ 1 \le r \le b,\ (r,b)=1$, it follows that
\begin{equation}
  E\left(0, \frac{r}{b}, 0\right) = \frac{1}{4} + \frac{i}{2} c_0\left(\frac{r}{b}\right),
\end{equation}

where $c_0(\frac{r}{b})$ is our cotangent sum from \eqref{eq:c0}.

\section{Central properties of the cotangent sum $c_0$}
Now we can state some crucial results concerning the cotangent sum $c_0$.
In \cite{Rassias}, M. Th. Rassias proved the following asymptotic formula:

\begin{theorem}\label{th:Rass}
  For $b \ge 2,\ b \in \mathbb{N}$, we have
  \begin{equation}
    c_0\left(\frac{1}{b}\right) = \frac{1}{\pi}b\log b - \frac{b}{\pi}(\log 2\pi - \gamma)
    + O(1).
  \end{equation}
\end{theorem}

Subsequently in \cite{Maier_Rassias}, M. Th. Rassias and H. Maier established an improvement,
or rather an asymptotic expansion, of Theorem \ref{th:Rass}.

\begin{theorem}\label{th:RassMa}
  Let $b,n \in \mathbb{N},\ b \ge 6N,$ with $N = \lfloor\frac{n}{2}\rfloor + 1.$
  There exist absolute real constants $A_1, A_2 \ge 1$ and absolute real
  constants $E_l,$ where $l \in \mathbb{N}$, with $|E_l| \le (A_1l)^{2l},$ such that for
  each $n \in \mathbb{N}$ we have
  \begin{equation}
    c_0\left(\frac{1}{b}\right) = \frac{1}{\pi}b\log b - \frac{b}{\pi}(\log2\pi - \gamma) -
    \frac{1}{\pi} + \sum_{l=1}^{n}E_l b^{-l} + R_{n}^{*}(b),
  \end{equation}
  where
  $$
    |R_{n}^{*}| \le (A_2n)^{4n}b^{-(n-1)}.
  $$
\end{theorem}

It is essential that both of these results were obtained using a common underlying idea
proposed in \cite{Rassias}. First of all, one can obtain the following
relation between sums of cotangents and sums with fractional parts:

\begin{equation}
  \mathop{\sum_{a \geq 1}}_{b\nmid a}
  \frac{b(1-2\{a/b\})}{a} =
  \mathop{\sum_{a \geq 1}}_{b\nmid a} \sum_{m = 1}^{b}
  \cot\left(\frac{\pi m}{b}\right)\frac{\sin(\frac{2\pi m}{b}a)}{a}.
\end{equation}

This relation provided the following proposition:

\begin{proposition}
  For every positive integer $b \geq 2$, we have
  \begin{equation}
    c_0\left(\frac{1}{b}\right) = \mathop{\sum_{a \geq 1}}_{b\nmid a}
    \frac{b(1-2\{a/b\})}{a}.
  \end{equation}
\end{proposition}

Then the difficulty lies in obtaining a good approximation of the sum $S(L;b)$ defined by

\begin{equation}
  S(L;b) := 2b \sum_{1\leq a\leq L}
  \frac{1}{a}\left\lfloor\frac{a}{b}\right\rfloor.
\end{equation}

The difference between the estimates from Theorems \ref{th:Rass} and \ref{th:RassMa}
is that stronger approximation techniques were applied in \cite{Maier_Rassias}
to obtain more information about $S(L;b)$. Namely, the generalized Euler
summation formula \eqref{eq:EulerSum} was used to improve the result of 
Theorem \ref{th:Rass}.

\begin{definition}
  If $f$  is a function that is differentiable at least $(2N+1)$ times in
  $[0, Z]$, let
  $$
    r_N(f, Z) = \frac{1}{(2N + 1)!}\int_{0}^{Z}(u -
    \lfloor u\rfloor + B)^{(2N + 1)}f^{(2N + 1)}(u)du,
  $$
  with the following notation:
  $$
    (u - \lfloor u\rfloor + B)^{(2N + 1)} =
    ((u - \lfloor u\rfloor) + B)^{(2N + 1)} :=
    \sum_{j=0}^{2N + 1}\binom{2N+1}{j}(u - \lfloor u\rfloor)^j B_{2N+1-j},
  $$
  where $B_{2j}$ are the Bernoulli numbers.
\end{definition}

\begin{theorem}[Generalized Euler summation formula]\label{th:Euler}
  Let $f$ be $(2N + 1)$ times differentiable in the interval $[0, Z]$. Then
  \begin{eqnarray}\label{eq:EulerSum}
    \lefteqn{
      \sum_{\nu = 0}^{Z} f(\nu) = \frac{f(0) + f(Z)}{2} + \int_{0}^{Z} f(u)du
    } \nonumber \\
    & & + \sum_{j=1}^{N}\frac{B_{2j}}{(2j)!}(f^{(2j - 1)}(Z) - f^{(2j - 1)}(0))
    + r_N(f, Z).
  \end{eqnarray}
\end{theorem}

Particularly, H. Maier and M. Th. Rassias used Theorem \ref{th:Euler} to obtain the following new
representation for $S(L;b)$.

\begin{theorem}
    For $N \in \mathbb{N}$, we have
    \begin{eqnarray}\label{eq:NewRepr}
      \lefteqn{
        S(L;b) = 2b \sum_{k \leq L/b} k\left(\log\frac{(k+1)b - 1}{kb - 1} +
        \frac{1}{2}F_1(k, b)\right)
      } \nonumber \\
      & & + 2b\sum_{j=1}^{N}\frac{B_{2j}}{2j}\sum_{k \leq L/b}kF_{2j}(k,b) +
      2br_N\left(f, \frac{L}{b}\right),
    \end{eqnarray}
    where the function $f$ satisfies
    \[ f(u) =
        \begin{cases}
          \frac{1}{u},\ if\ u\geq1\\
          0,\ if\ u=0
        \end{cases}
    \]
    and $f\in C^{\infty}([0, \infty))$ with $f^{(j)}(0) = 0$ for $j\leq 2N + 1$.
\end{theorem}

The new form of $S(L;b)$ from \eqref{eq:NewRepr} lead essentially to the proof of Theorem
\ref{th:RassMa}.

Furthermore, H. Maier and M. Th. Rassias obtained even more interesting results concerning
$c_0(\frac{r}{b})$ for a fixed arbitrary positive integer value of $r$ and for
large integer values of $b$, which give us a deeper understanding of our cotangent
sum for almost all values of $r$ and $b$.

\begin{proposition}\label{prop:c_0gen}
  For $r,b \in \mathbb{N}$ with $(r,b)=1$, it holds that
  \begin{equation}
    c_0\left(\frac{r}{b}\right) = \frac{1}{r}c_0\left(\frac{1}{b}\right) - \frac{1}{r}Q\left(\frac{r}{b}\right),
  \end{equation}
  where
  $$
    Q\left(\frac{r}{b}\right) = \sum_{m=1}^{b-1} \cot\left(\frac{\pi mr}{b}\right)
    \left\lfloor\frac{rm}{b}\right\rfloor.
  $$
\end{proposition}

\begin{theorem} \label{th:RassMaier}
  Let $r,b_0 \in \mathbb{N}$ be fixed, with $(b_0, r) = 1$. Let $b$ denote a
  positive integer with $b \equiv b_0 (mod\ r)$. Then, there exists a constant
  $C_1 = C_1(r, b_0)$, with $C_1(1, b_0) = 0$, such that
  \begin{equation}
    c_0\left(\frac{r}{b}\right) = \frac{1}{\pi r}b\log b - \frac{b}{\pi r}
    (\log2\pi - \gamma) + C_1b + O(1)
  \end{equation}
  for large integer values of $b$.
\end{theorem}

The function $c_0$ is also thoroughly studied in the papers of S. Bettin and J. Conrey
\cite{Bettin_Conrey_2011, Bettin_Conrey_2013}, where they have established a reciprocity formula for it. However, before we state the formula itself, we must give several definitions.

For $a \in \mathbb{C}$ and Im $(s) > 0$, consider
$$
  \mathscr{S}_a(s) := \sum_{n=1}^{\infty} \sigma_a(n)e(ns)e^{2\pi ins},
$$

$$
  E_a(s) := 1 + \frac{2}{\zeta(-a)}\mathscr{S}_a(s).
$$
It is worth mentioning that for $a = 2k + 1$, $k \in \mathbb{Z}_{\geq 1}$,
$E_a$ is the well known Eisenstein series of weight $2k+2$.
\begin{definition}
  For $a \in \mathbb{C}$ and $Im(s) > 0$, define the function
  \begin{equation}
    \psi_a(s) := E_{a+1}(s) - \frac{1}{s^{a+1}}E_{a+1}\left(-\frac{1}{s}\right).
  \end{equation}
\end{definition}

For $a = 2k$, $k \geq 2$, the function $\psi_a(s)$ is equal to zero, because of the
modularity property of the Eisenstein series. Unfortunately, it is not true for
other values of $a$, but the functions $\psi_a(s)$ have some remarkable properties,
which were described in detail by S. Bettin and J. Conrey.

Now we can state the theorem proven in paper \cite{Bettin_Conrey_2011}.
\begin{theorem}\label{th:recip}
  The function $c_0$ satisfies the following reciprocity formula:
  \begin{equation}
    c_0\left(\frac{r}{b}\right) + \frac{b}{r}c_0\left(\frac{b}{r}\right) - \frac{1}{2\pi r} =
    \frac{i}{2}\psi_0\left(\frac{r}{b}\right).
  \end{equation}
\end{theorem}

This result implies that the value of $c_0(\frac{r}{b})$ can be computed  within a
prescribed accuracy in a polynomial of $\log b$.

S. Bettin and J. Conrey highlighted that the reciprocity formula from \ref{th:recip} is
very similar to that of the Dedekind sum \ref{th:Dedekind_recip}. We will
consider Dedekind sums in more detail in section \ref{sec:Dedekind} of this
paper.

In \cite{Bettin_Conrey_2013} the result for $c_0$ was generalized for the sums

\begin{equation}\label{eq:cot_general}
  c_a\left(\frac{r}{b}\right) :=
  b^a \sum_{m=1}^{b-1}\cot\left(\frac{\pi mr}{b}\right)\zeta\left(-a, \frac{m}{b}\right),
\end{equation}
where $\zeta(s, x)$ is the Hurwitz zeta function.

\subsection{Ellipse}

\begin{figure}[ht]
\centering
  \includegraphics[scale=1.0]{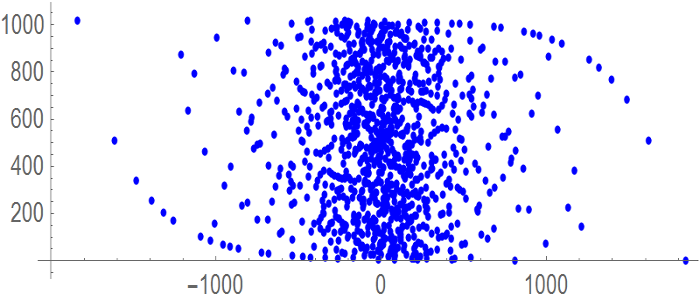}
  \caption{axis $Ox: r$, axis $Oy: c_0$, $b = 1021$}
  \label{fig:HD1}
\end{figure}

\begin{figure}[ht]
\centering
  \includegraphics[scale=1.0]{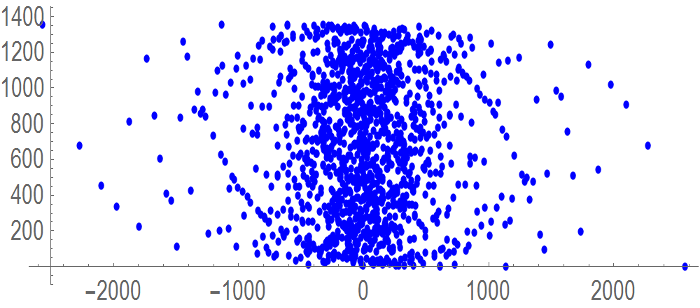}
  \caption{axis $Ox: r$, axis $Oy: c_0$, $b = 1357$}
  \label{fig:HD2}
\end{figure}

It is interesting to mention that if one examines the graph of $c_0\left(\frac{r}{b}\right)$ for hundreds of integer values of $b$ by the use of MATLAB, the resulting figures \ref{fig:HD1} and \ref{fig:HD2} always have a shape similar to an ellipse.

In \cite{Maier_Rassias} H. Maier and M. Th. Rassias tried to explain this phenomenon and obtained an important result, which establishes the equidistribution of certain normalized cotangent sums with respect to a positive measure. This is presented in the following theorem.

\begin{definition}\label{defg}
  For $z \in \mathbb{R}$, let
  
  $$
  F(z)= meas\{ x \in [0,1] : g(x)\leq z\},
  $$
  where “meas” denotes the Lebesgue measure,
  
  $$
    g(x):= \sum_{l=1}^{+\infty} \frac{1-2\{lx\}}{l},
  $$
  and 
  $$
  C_0(\mathbb{R}) = \{f \in C(\mathbb{R}) : \forall \varepsilon > 0,
  \exists ~ \text{a compact set} ~ \mathcal{K} \subset \mathbb{R},
\text{such that} ~ |f(x)| < \varepsilon, \forall x \notin \mathcal{K}\}.
  $$
\end{definition}
\textbf{Remark.} The convergence of the above series has been investigated by R.Bret\`eche and G.Tenenbaum
(see Theorem \ref{th:Breteche-Tenenbaum}). It depends on the partial fraction expansion of the number $x$.
\begin{theorem} \label{th:Ellipses}
  \textnormal{(i)} $F$ is a continuous function of $z$.
  
  \textnormal{(ii)}  Let $A_0$, $A_1$ be fixed constants, such that $1/2 < A_0 < A_1 < 1$. Let also
  $$
  H_{k}=\int_0 ^1 \left(\frac{g(x)}{\pi}\right)^{2k} dx,
  $$
  so $H_k$ is a positive constant depending only on $k$, $k \in \mathbb{N}$.
  
  There is a unique positive measure $\mu$ on $\mathbb{R}$ with the following properties:
  
  \textnormal{(a)} For $\alpha < \beta \in \mathbb{R}$ we have
  
  $$
  \mu([\alpha, \beta]) = (A_1 − A_0)(F(\beta) − F(\alpha)).
  $$
  
  \textnormal{(b)}
  
  $$
  \int x^k d\mu = \begin{cases}
  (A_1 − A_0)H_{k/2},  & \mbox{for even } k \\
  0, & \mbox{otherwise.}
\end{cases}
  $$
  
  \textnormal{(c)}
  For all $f \in C_0(\mathbb{R})$, we have
  
  $$
  \lim \limits_{b\rightarrow+\infty} \frac{1}{\phi(b)} \sum \limits_{\substack{r:(r,b)=1,\\ A_0 b\leq r \leq A_1 b}}
  f \left(\frac{1}{b} c_0\left(\frac{r}{b}\right)\right) =
  \int f d\mu,
  $$
  where $\phi(\cdot)$ denotes the Euler phi-function.
  
\end{theorem}

\textbf{Outline of the proof.}

In \cite{Maier_Rassias} H. Maier and M. Th. Rassias proved Theorem \ref{th:RassMa}, which constitutes an improvement of their earlier Theorem \ref{th:Rass}.

Additionally, they investigated the cotangent sum $c_0\left(\frac{r}{b}\right)$ for a fixed arbitrary positive
integer value of $r$ and for large integer values of $b$ and proved Theorem \ref{th:RassMaier} as well as the following results

\begin{theorem} \label{th:Moments}
   Let $k \in \mathbb{N}$ be fixed. Let also $A_0, A_1$ be fixed constants such that
   
  $\frac{1}{2} < A_0 < A_1 < 1$. Then there exists a constant $E_k > 0$,
  depending only on $k$, such that
  
   \textnormal{(a)}
   
  $$
    \mathop{\sum_{r:(r,b)=1}}_{A_0b≤r≤A_1b}  Q\left(\frac{r}{b}\right)^{2k} = E_k
    (A^{2k+1}_1 − A^{2k+1}_0 )b^{4k}\phi(b)(1 + o(1))
    \;\;\;\; (b \rightarrow +\infty),
  $$
  
  \textnormal{(b)}
  
  $$
    \mathop{\sum_{r:(r,b)=1}}_{A_0b≤r≤A_1b}  Q\left(\frac{r}{b}\right)^{2k-1} =o(b^{4k-2}\phi(b))
    \;\;\;\; (b \rightarrow +\infty),
  $$
 
 \textnormal{(c)}
 
  $$
    \mathop{\sum_{r:(r,b)=1}}_{A_0b≤r≤A_1b} c_0\left(\frac{r}{b}\right)^{2k} = H_k
    (A_1 - A_0)b^{2k}\phi(b)(1 + o(1))
    \;\;\;\; (b \rightarrow +\infty),
  $$
  
  \textnormal{(d)}
  
   $$
    \mathop{\sum_{r:(r,b)=1}}_{A_0b≤r≤A_1b} c_0\left(\frac{r}{b}\right)^{2k-1} = o(b^{2k-1}\phi(b))
    \;\;\;\; (b \rightarrow +\infty),
  $$
   with
  $$
    E_k=\frac{H_{k}}{(2k+1)}.
  $$
 
\end{theorem}

Using the method of moments, one can deduce detailed information about the distribution of the values of $c_0\left(\frac{r}{b}\right)$, where $A_0 b \leq r \leq A_1 b$ and $b \rightarrow + \infty$. Namely, one can prove Theorem \ref{th:Ellipses}.

\section{The maximum of $c_0$ in rational numbers in short intervals}

In this section we consider some results about the maximum of $c_0$ in rational
numbers in short intervals. More precicely, consider the following definition:

\begin{definition}\label{def:M1}
  Let $0 < A_0 < 1$, $0 < C < 1/2$. For $b \in \mathbb{N}$ we set
  $$
    ∆ := ∆(b, C) = b^{-C}.
  $$

  We define
  $$
    M(b, C, A_0) := \mathop{\text{max}}_{A_0b≤r<(A_0+∆)b} \left|c_0\left(\frac{r}{b}\right)\right|.
  $$
\end{definition}

In \cite{maximum} H. Maier and M. Th. Rassias proved the following
theorems.

\begin{theorem}\label{th:M1}
  Let $D$ satisfy $0 < D < \frac{1}{2} − C$.
  Then we have for sufficiently large $b$:
  $$
    M(b, C, A_0) \geq\frac{D}{\pi}b\log b.
  $$
\end{theorem}
An important part of the proof is the following key proposition.

\begin{proposition} \label{prop: M1_1}
  Let $\langle a_0; a_1, a_2, \ldots, a_n\rangle$ be the continued fraction
  expansion of $\frac{\bar{r}}{b} \in \mathbb{Q}$. Moreover, let $\frac{u_l}{v_l}$
  be the $l$-th partial quotient of $\frac{\bar{r}}{b}$. Then
  $$
    c_0\left(\frac{r}{b}\right)=-b\sum_{1 \leq l\leq n}\frac{(-1)^l}{v_l}
    \left(\left( \frac{1}{\pi v_l}\right) +
    \psi\left(\frac{v_{l-1}}{v_l} \right)\right).
  $$

  Here $\psi$ is an analytic function satisfying

  $$
    \psi(x) = - \frac{\log(2\pi x)-\gamma}{\pi x} + O(\log x),
    \; (x \rightarrow 0).
  $$
\end{proposition}

The proposition was proven in \cite{Bettin_Moments} by S. Bettin.

\begin{definition}
  Let $\Delta$ be as in Definition \ref{def:M1} and $\Omega > 0$. We set

$$
  N(b, \Delta, \Omega) := \#\{r : A_0 b \leq r < (A_0 + \Delta)b,
  \; |\bar{r}| ≤ \Omega b\}.
$$

\end{definition}

Another key proposition, proven in \cite{maximum}, is the following:

\begin{proposition} \label{prop: M1_2}
  Let $\varepsilon <0$ be such that
  $$
    D + \varepsilon < \frac{1}{2} - C.
  $$
  Set
  $$
    \Omega := b^{-(D+\varepsilon)}.
  $$
  Then for sufficiently large $b$ it holds 
  $$
    N(b, \Delta, \Omega)>0.
  $$
\end{proposition}

Let $\{\frac{u_i}{v_i}\}_{i=1}^s$ be the sequence of partial fractions of
such $\frac{\bar{r}}{b}$. From

$$
  \Omega\geq\frac{\bar{r}}{b}\geq\frac{1}{v_1+1}
$$
we obtain

$$v_1+1\geq \Omega^{-1}.$$

Then by proposition \ref{prop: M1_2} we have

$$
  \sum_{l>1}  \left(\left( \frac{1}{\pi v_l}\right) +
  \psi\left(\frac{v_{l-1}}{v_l} \right)\right) <
  2\varepsilon \log b, \; \text{for}  \; b\geq b_0(\varepsilon).
$$
Therefore,

$$
  \left| c_0\left(\frac{r}{b}\right)\right|\geq \frac{2}{\pi}\log(\Omega^{-1}(1+o(1)))
  \; \; b\rightarrow+\infty.
$$

This proves Theorem \ref{th:M1}.

\begin{theorem} \label{th:M2}
  Let C be as in Theorem \ref{th:M1} and let $D$ satisfy $0 < D < 2 − C - E$,
  where $E \geq 0$ is a fixed constant. Let $B$ be sufficiently large.
  Then we have
  $$
    M(b, C, A_0) \leq\frac{D}{\pi}b\log b
  $$
  for all $b$ with $B \leq b < 2B$, with at most $B^E$ exceptions.
\end{theorem}
\textbf{Proof.}

We will need the following proposition.

\begin{proposition} \label{prop:M2_1}
  Let $ \varepsilon> 0, \; B \geq B\varepsilon, \; B < b \leq 2B$.
  For $1 \leq r < b, \; (r, b) = 1$, let $\{\frac{u_i}{v_i}\}_{i=1}^s$ be the
  sequence of partial fractions of $\frac{\bar{r}}{b}$. Then there are at most $3$
  values of $l$ for which
  $$
    \frac{1}{v_l}\psi\left(\frac{v_{l-1}}{v_l}\right)\geq \log\log b,
  $$
  and at most one value of $l$ for which

  $$
    \frac{1}{v_l}\psi\left(\frac{v_{l-1}}{v_l}\right)\geq \varepsilon\log b.
  $$
\end{proposition}
\textbf{Proof.} Let $l_i (i = 1, 2, 3, 4)$ be such that

$$
  \frac{1}{v_l}\psi\left(\frac{v_{l-1}}{v_l}\right)\geq \log\log b.
$$
Then we have

$$
  v_{l_1} \geq \log\log b,\ \ v_{l_2} \geq \exp(v_{l_1}) \geq \log b,
$$
$$
  v_{l_3} \geq \exp(v_{l_2}) \geq b,\ \ v_{l_4} \geq \exp(v_{l_3}) \geq \exp(b),
$$
in contradiction to $v_s \leq b$.
\bigskip

In the same manner we obtain from $v_{l_j} \geq \varepsilon \log b, j = 1, 2 :$
$$v_s \geq \exp(\exp((\log b)^{\varepsilon})) > b .$$

\bigskip

Assume $\varepsilon > 0$ to be fixed but arbitrarily small, $Z > 0$ fixed but
arbitrarily large.

\begin{definition}
  By Proposition \ref{prop:M2_1} there is at most one value of $l$ for which
  $$
    \frac{1}{v_l}\psi\left(\frac{v_{l-1}}{v_l}\right)\geq \varepsilon\log b.
  $$
  In case of the existence of $l$, we write

  $$
    u_{l−1}(r, b) = u_{l−1}
  $$
  and

  $$
    v_{l−1}(r, b) = v_{l−1}.
  $$
  Then for $s, t$ with $1 \leq s, t ≤ Z, (s, t) = 1$, and for fixed $\theta$ with
  $0 < \theta < 1$,

  $$
    \mathcal{F}(s, t): = \left\{ (b,r,\bar{r}): B \leq b < 2B, \; A_0 b \leq
    r ≤ (A_0 + \Delta)b, \; \left| \frac{\bar{r}}{b}-\frac{s}{t}\right|\leq
    \theta, \; r\bar{r} \equiv 1 (\mathrm{mod}\ b) \right\}.
  $$
\end{definition}

Now we can formulate a proposition.

\begin{proposition} \label{prop:M2_2}
$$
  \sum _{B \leq b < 2B} N(b, \Delta, \Omega) \leq \sum_{1\leq s,t \leq Z}
  \left|\mathcal{F}(s, t)\right|
$$
\end{proposition}

By Dirichlet’s approximation theorem there is $(C_0, D_0) \in \mathbb{Z}^2$ with
$1 \leq D_0 \leq B^2$, $(C_0, D_0) = 1,$ such that
\begin{equation} \label{eq:M2_1}
  \left| A_0^{-1} - \frac{C_0}{D_0}\right| \leq \frac{1}{D_0 B^2}.
\end{equation}

Let us estimate the cardinality of the set $\mathcal{F}(s, t)$.

From $A_0 b \leq r ≤ (A_0 + \Delta)b$, the definition of $\bar{r}$, and
\eqref{eq:M2_1}, we obtain

\begin{equation}
  r\bar{r}= y\left(\frac{C_0}{D_0}r + \frac{u}{D_0}\right),
  \text{with} \; y \in \mathbb{Z},
\end{equation}
which after multiplication by $C_0 D_0$ becomes

\begin{equation} \label{eq:M2_2}
  (C_0 y − D_0\bar{r})(C_0 r − u) = −D_0(C_0 − \bar{r}u ).
\end{equation}

If $C_0 −\bar{r}u \neq 0$, one can deduce from the well-known estimate for the
number of divisors of an integer that for a given pair $(\bar{r}, u)$, there are
at most $O(B^{\varepsilon})$ pairs $(r, y)$ such that \eqref{eq:M2_2} holds.

There are at most $O(B^{\varepsilon})$ pairs $(\bar{r}, u)$ such that
$C_0 − \bar{r}u = 0$. Thus we obtain

\begin{equation} \label{eq:M2_3}
  \left|\mathcal{F}(s, t)\right| = O(B^{2+2\varepsilon - C}\theta).
\end{equation}

From Proposition \ref{prop:M2_2} and \eqref{eq:M2_3} we obtain for
$\Omega=\theta B$

\begin{equation} \label{eq:M2_4}
  \sum _{B \leq b < 2B} N(b, \Delta, \Omega) = O(B^{2+2\varepsilon - C}\theta).
\end{equation}

We now apply \eqref{eq:M2_4} with $\theta = B−D_0$, where

$$
  D > D' > 2 − C − E.
$$

If we choose $\varepsilon>0$ sufficiently small, then we  conclude from
\eqref{eq:M2_4} the following:

For all $b$ with $B \leq b < 2B$ we have, with at most $B^E$ exceptions:

$$
  N(b, B^{1-C}, \Omega) =0.
$$
Thus,

\begin{equation}
  \frac{1}{v_l}\psi\left(\frac{v_{l-1}}{v_l}\right) \leq
  \frac{D'}{\pi v_{l-1}} b \log b (1+ o(1)),\
  \forall l \leq Z.
\end{equation}

The result of Theorem \ref{th:M2} follows now from Propositions \ref{prop: M1_1}
and \ref{prop:M2_1}.

\section{The function $g(x)$ and moments of $c_0$}

There is an interesting connection between the cotangent sums $c_0$ and the function
\begin{equation}\label{eq:g(x)}
    g(x):= \sum_{l=1}^{+\infty} \frac{1-2\{lx\}}{l},
\end{equation}
which, as we mentioned above, naturally appeared in the investigation of
the moments of $c_0$ and of the sum $Q\left(\frac{r}{b}\right)$,
which is related to $c_0$ by Proposition \ref{prop:c_0gen}. To be precise, this series is related to $c_0$ by the important Theorem \ref{th:Moments}.

Later S. Bettin in his paper \cite{Bettin_Moments} extended the result of Theorem \ref{th:Moments}
and proved the following:

\begin{theorem}
    Let $b \geq 1$ and $k \geq 0$. Then
    \begin{equation}
      \frac{1}{\phi(b)}\mathop{\sum^{b}_{r=1}}_{(r, b)=1}c_0\left(\frac{r}{b}\right)^k =
      H_kb^k + O_{\varepsilon}(b^{k-1+\varepsilon}(Ak\log b)^{2k}),
    \end{equation}
    for some absolute constant $A>0$ and any $\varepsilon > 0$.

    Moreover, if $0 \leq A_0 < A_1 \leq 1$, then we have
    \begin{equation}
      \frac{1}{\phi(b)}\mathop{\sum_{(r, b)=1}}_{A_0 < \frac{r}{b} < A_1}
      c_0\left(\frac{r}{b}\right)^k = (A_1 - A_0)H_kb^k +
      O_{\varepsilon}(b^{k-\frac{1}{2}+\varepsilon}(Ak\log b)^{2k}).
    \end{equation}
\end{theorem}

The function $g(x)$ is interesting not only in connection to the study of the
cotangent sums $c_0$, but also in its own right. For example, it is also studied in
\cite{Tenenbaum} by R. Bret\`eche and G. Tenenbaum.

\begin{theorem} \label{th:Breteche-Tenenbaum}
  For each $x \in \mathbb{Q}$ the series $g(x)$ converges.

  For $x \in \mathbb{R}\setminus\mathbb{Q}$, the series $g(x)$ converges if and
  only if the series
  $$
    \sum_{m\geq 1}(-1)^m \frac{\log q_{m+1}}{q_m}
  $$
  converges, where $(q_m)_{m \geq 1}$ denotes the sequence of partial
  denominators of the continued fraction expansion of $x$.
\end{theorem}

\textbf{Proof.} The statement of the theorem is part of Theorem $4.4$ of
the paper by R. Bret\`eche and G. Tenenbaum in \cite{Tenenbaum}.
\bigskip

The function $g(x)$ also has the following property:
\begin{theorem}
  The series
  $$
    g(x) = \sum_{l=1}^{+\infty} \frac{1-2\{lx\}}{l}
  $$
  converges almost everywhere in $[0, 1)$.
\end{theorem}

The function $g(x)$ was also of interest to L. B\'aez-Duatre, M. Balazard, B. Landreau and
E. Saias. In \cite{Balazard&Co} they studied the function

$$
  A(\lambda):=\int_0^{+\infty} \{t\}\{\lambda t\} \frac{dt}{t^2}.
$$

and proved the following theorem.

\begin{theorem}
  Let $\lambda > 0$ be such that the series $g(\lambda)$ converges. Then the
  series $g(\frac{1}{\lambda})$ converges too, and we have:

  $$
    A(\lambda)=\frac{1-\lambda}{2}\log \lambda + \frac{\lambda+1}{2}(\log 2\pi -
    \gamma) - g(\lambda) - \lambda g\left(\frac{1}{\lambda}\right).
  $$
\end{theorem}

\bigskip

Now let us show an important property of $g(x)$, which was proven in
\cite{moments}.

\begin{theorem} \label{th: moments_g(x)}
  There are constants $c_1, c_2 > 0$, such that
  $$
    c_1 \Gamma(2k + 1) \leq \int_0 ^1 g(x)^{2k} dx \leq c_2 \Gamma(2k + 1)
  $$
  for all $k \in \mathbb{N}$, where $\Gamma(\cdot)$
  stands for the gamma function.
\end{theorem}
\textbf{Sketch of a proof.} Let us consider the continued fraction expansion of $x$

$$
x = [a_0(x); a_1(x),\dots,a_k(x),\dots].
$$

The $a_k(x)$ are obtained via the Gauss map $\alpha$, defined by

$$ 
\alpha(x)=\left\{\frac{1}{x}\right\}, \alpha_k(x)=\alpha(\alpha_{k-1}(x)), a_k(x)=\left\lfloor \frac{1}{\alpha_{k-1}(x)} \right\rfloor.
$$

\begin{definition}
   Let $x \in X = (0, 1) \backslash \mathbb{Q}$. Let also
   
   $$
   \beta_k(x)=\alpha_0(x)\alpha_1(x)\dots\alpha_{k}(x), \; \beta_{-1}(x)=1,
   $$
   $$
   \gamma_k(x)=\beta_{k-1}(x)\log\frac{1}{\alpha_k(x)}, \text{where} \; k \geq 0,
   $$
   
   so that $\gamma_0(x)=\log\frac{1}{x}$.
   
   The number $x$ is called a Wilton number if the series
   
   $$
   \sum_{k \geq 0} (-1)^{k}\gamma_k(x)
   $$
   
   converges.
   
   Wilton’s function $\mathcal{W}$ is defined by
   
   $$
   \mathcal{W} = \sum_{k \geq 0} (-1)^{k}\gamma_k(x)
   $$
   
   for each Wilton number $x$.
\end{definition}

M. Balazard and B. Martin proved in \cite{Balazard-Martin} the following proposition:

\begin{proposition}
  There is a bounded function $H : (0, 1) \rightarrow \mathbb{R}$, which is continuous at every irrational number, such that

\begin{equation} \label{eq:W(x)_0}
    g(x)= \mathcal{W}(x)+H(x)
\end{equation}

almost everywhere. Also a number $x \in X$ is a Wilton number if and only if $\alpha(x)$ is a Wilton number. In this case, we have:

\begin{equation} \label{eq:W(x)_1}
\mathcal{W}(x)=l(x) - T\mathcal{W}(x),
\end{equation}

where

$$
l(x)=\log\frac{1}{x}
$$
and the operator T is defined by

$$
Tf(x)= x f(\alpha(x)).
$$
\end{proposition}

One can express (\ref{eq:W(x)_1}) as

\begin{equation} \label{eq:W(x)_2}
    l(x)=(1+T)\mathcal{W}(x).
\end{equation}

The main idea in the evaluation of 

$$
\int_0 ^1 g(x)^{2k} dx
$$
is to solve the operator equation (\ref{eq:W(x)_2}) for $\mathcal{W}(x)$, which is:

\begin{equation} \label{eq:W(x)_3}
    \mathcal{W}(x)=(1+T)^{-1}l(x).
\end{equation}

An idea which has long been used in functional analysis for the case when $T$ is a differential operator is to express the right-hand side of (\ref{eq:W(x)_3}) as a Neumann series, which is obtained by the geometric series identity, i.e.

$$
(1+T)^{-1}=\sum_{k=0}^{+\infty}(-1)^k T^k.
$$

Thus one can approximate $\mathcal{W}(x)$ by
\begin{equation}
    \mathcal{L}(x,n)=\sum_{k=0}^{n} (-1)^k (T^k l)(x).
\end{equation}

\begin{definition}
  The measure $m$ is defined by
  
  $$
  m(\mathcal{E}) = \frac{1}{\log 2} \int \limits_{\mathcal{E}} \frac{dx}{1+x},
  $$
  where $\mathcal{E}$ is any measurable subset of $(0, 1)$.
 
\end{definition}

\begin{proposition} \label{prop:MMY}
  For $f \in L^p$ we have
  $$
  \int\limits_0^1 |T^n f(x)|^p dm(x) \leq g^{(n-1)p} \int\limits_0^1 |f(x)|^p dm(x),
  $$
  where 
  
 $$
 g=\frac{\sqrt{5}-1}{2}<1.
 $$
\end{proposition}
\textbf{Proof.} Marmi, Moussa, and Yoccoz, in their paper \cite{MMY}, consider a generalized continued fraction algorithm,
depending on a parameter $\alpha$, which becomes the usual continued fraction algorithm for the choice $\alpha = 1$. The
operator $T_v$ is defined in $(2.5)$ of \cite{MMY} and becomes $T$ for $\alpha = 1, v=1$.
Then, proposition \ref{prop:MMY} is the content of formulas
$(2.14)$, $(2.15)$ of \cite{MMY}.

\bigskip

Using standard techniques of functional analysis one can prove that

$$
    \lim \limits_{n \rightarrow \infty} \int \limits_{0}^{1/2}\left| \mathcal{L}(x,n)^L-\mathcal{W}(x)^L\right| dx =0.
$$

One can eventually prove that

$$
\int \limits_{0}^{1/2} \mathcal{L}(x,n)^{2k} dx = \left(\int \limits_{0}^{1/2} l(x)^{2k}\right) (1+o(1)) = \Gamma(2k+1) (1+o(1)), \; (k \rightarrow +\infty).
$$

The order of magnitude of
$$
\int \limits_{0}^{1/2} g(x)^{2k} dx
$$

now follows from (\ref{eq:W(x)_0}), by the binomial theorem, since $H(x)$ is a bounded function.

\bigskip

\begin{corollary}
  The series
  $$
    \sum_{k\geq 0} \frac{H_{k}}{(2k)!}x^k
  $$
  has radius of convergence $\pi^2$.
\end{corollary}

H. Maier and M. Th. Rassias proved in \cite{Asymptotic_moments_1} an improvement of Theorem \ref{th: moments_g(x)}, by establishing an asymptotic  result for the corresponding integral. Namely, they proved the following theorem.

\begin{theorem}
  Let
  $$
  A= \int\limits_0^{\infty} \frac{\{t\}^2}{t^2}~ dt
  $$
  
  and $K \in \mathbb{N}$. There is an absolute constant $C > 0$ such that
  $$
  \int \limits_{0}^{1} |g(x)|^{K} dx = 2e^{-A} \Gamma(K+1)(1+O(\exp(-CK)))
  $$
  
  for $K \rightarrow \infty$.
\end{theorem}

In \cite{Asymptotic_moments_2}, they improved this result settling also the general case of arbitrary exponents $K$.

\begin{theorem}
  Let $K \in \mathbb{R}, ~~ K>0$. There is an absolute constant $C > 0$ such that
  $$
  \int \limits_{0}^{1} |g(x)|^{K} dx = \frac{e^{\gamma}}{\pi} \Gamma(K+1)(1+O(\exp(-CK)))
  $$
  
  for $K \rightarrow \infty$, where $\gamma$ is the Euler-Mascheroni constant.
\end{theorem}

\section{Dedekind sums}\label{sec:Dedekind}

Dedekind sums have applications in many fields of mathematics, especially in
number theory.  These sums appear in R. Dedekind's study of the
function
\begin{equation}
  \eta(s) = e^{\frac{\pi is}{12}}\prod_{m=1}^{\infty}(1 - e^{2\pi ims}),
\end{equation}
where Im $s > 0$.

\begin{definition}
  Let $r,b$ be integers, $(r, b) = 1$, $k \geq 1$. Then the Dedekind sum $s(\frac{r}{b})$ is defined as follows
  \begin{equation}
    s\left(\frac{r}{b}\right) := \sum_{\mu = 1}^{b}\left(\left(\frac{r\mu}{b}\right)\right)\left(\left(\frac{\mu}{b}\right)\right),
  \end{equation}
  where $((\cdot))$ is the sawtooth function defined as follows:
  \begin{equation}
    ((x)) :=
      \begin{cases}
        x - [x] - \frac{1}{2} & \quad \text{if } $x$ \text{ is not an integer,}\\
        0 & \quad \text{if } $x$ \text{ is an integer.}\\
      \end{cases}
  \end{equation}
\end{definition}

It is a fascinating fact that the Dedekind sum can also be expressed as a sum of
cotangent products:
\begin{proposition}\label{prop:equiv_form}
  $$
    s\left(\frac{r}{b}\right) = -\frac{1}{4b}\sum_{m=1}^{b-1}\cot\left(\frac{\pi m}{b}\right)
    \cot\left(\frac{\pi mr}{b}\right).
  $$
\end{proposition}

It is a well-known fact that Dedekind sums satisfy a reciprocity formula:
\begin{theorem}[Dedekind sums' reciprocity formula]\label{th:Dedekind_recip}
  \begin{equation}
    s\left(\frac{r}{b}\right) + s\left(\frac{b}{r}\right) - \frac{1}{12rb} =
    \frac{1}{12}\left(\frac{r}{b} + \frac{b}{r} -3\right).
  \end{equation}
\end{theorem}

We shall not present the proof of Proposition \ref{prop:equiv_form} and of Theorem
\ref{th:Dedekind_recip}. The interested reader can find these proofs, as well as more fundamental facts
concerning Dedekind sums in the now famous book by Rademacher and
Grosswald \cite{Dedekind_Sums}.

It is interesting to study the relation between the cotangent sums $c_0$ and the
Dedekind sums. We've already considered $c_a$ for arbitrary $a \in \mathbb{C}$
(see (\ref{eq:cot_general})). S. Bettin in \cite{Bettin_Rademacher} proved the very interesting result 
that $c_{-1}$ is a Dedekind sum up to a constant:

\begin{proposition} It holds that
  \begin{equation}
    s\left(\frac{r}{b}\right) = \frac{1}{2\pi}c_{-1}\left(\frac{r}{b}\right).
  \end{equation}
\end{proposition}

\section{Sums appearing in the Nyman-Beurling criterion for the Riemann
Hypothesis containing the M\"{o}bius function}\label{sec:SumsWithMo}

In a recent paper \cite{Sums_Estimates} H. Maier and M. Th. Rassias investigated 
the following sums  \eqref{eq:g(x)}
\begin{equation}\label{eq:sums_with_mu}
  \sum_{n \in I} \mu(n)g\left(\frac{n}{b}\right),
\end{equation}
for a suitable interval $I$, where $\mu(n)$ is the M\"{o}bius function and the function $g(x)$ is as defined in Definition \ref{defg}.

These sums appear in the study of the integral
\begin{equation}\label{eq:last_int}
  \int_{-\infty}^{+\infty}
  \left|\zeta\left(\frac{1}{2} + it\right)\right|^2
  \left|D_N\left(\frac{1}{2} + it\right)\right|^2
  \frac{dt}{\frac{1}{4}+t^2}
\end{equation}

Particularly we could express \eqref{eq:last_int} using formulas
\eqref{eq:V_N} and \eqref{eq:exp_integral} and Proposition \ref{prop:int_to_vas}, as follows:
\begin{align*}
  \int_{-\infty}^{+\infty} &
  \left|\zeta\left(\frac{1}{2} + it\right)\right|^2
  \left|D_N\left(\frac{1}{2} + it\right)\right|^2
  \frac{dt}{\frac{1}{4}+t^2} \\
  & = \sum_{1 \leq r,b \leq N} \mu(r)\mu(b)
  \left(1 - \frac{\log r}{\log N}\right)
  \left(1 - \frac{\log b}{\log N}\right) \\
  & \times \left[\frac{\log2\pi - \gamma}{2}\left(\frac{1}{r} +
  \frac{1}{b}\right) + \frac{b-r}{2rb}\log\frac{r}{b} -
  \frac{\pi}{2rb}\left(V\left(\frac{r}{b}\right) +
  V\left(\frac{b}{r}\right)\right)\right]
\end{align*}

If we expand the last equation, we will obtain the following sum for fixed $b$
$$
  \sum_{n \in I} \mu(n)\left(1 - \frac{\log n}{\log N}\right)
  \frac{1}{n}V\left(\frac{n}{b}\right),
$$
which is equal to
$$
  \sum_{n \in I} \mu(n)\left(1 - \frac{\log n}{\log N}\right)
  g\left(\frac{n}{b}\right).
$$

In \cite{Sums_Estimates} H. Maier and M. Th. Rassias proved the following result concerning the
sums \eqref{eq:sums_with_mu}:

\begin{theorem}
Let $0 \leq \delta \leq D/2$, $b^{2\delta} \leq B \leq b^D$, where
$b^{-\delta}\leq \eta \leq 1$. Then there is a positive constant $\beta$
depending only on $\delta$ and $D$, such that
\begin{equation}
  \sum_{Bb\leq n \leq (1 + \eta)Bb} \mu(n)g\left(\frac{n}{b}\right)
  = O((\eta Bb)^{1-\beta}).
\end{equation}
\end{theorem}

Finally, the above result was recently improved by the same authors in \cite{maier_rassias_jfa}, by proving the following theorem:

\begin{theorem}
Let $D\geq 2$. Let $C$ be the number which is uniquely determined by
$$C\geq \frac{\sqrt{5}+1}{2},\ 2C-\log C - 1- 2\log 2=\frac{1}{2}\log 2\:.$$
Let $v_0$ be determined by
$$v_0\left(1-\left(1+2\log 2\left(C+\frac{\log 2}{2}\right)^{-1}\right)^{-1}+2+\frac{4}{\log 2} C\right)=2\:.$$
Let $z_0:=2-\left(2+\frac{4}{\log 2} C\right) v_0$. Then for all $\varepsilon >0$ we have
$$\sum_{b^D\leq n< 2b^D} \mu(n)g\left(\frac{n}{b}\right)\ll_{\varepsilon} b^{D-z_0+\varepsilon}\:.$$
\end{theorem}

\end{document}